\documentclass{amsart}
\usepackage{amsmath}
\usepackage{amssymb}
\usepackage{amsfonts}

\newtheorem{theorem}{Theorem}
\theoremstyle{plain}

\newtheorem{lemma}{Lemma}

\newtheorem{proposition}{Proposition}

\numberwithin{equation}{section}
\input{tcilatex}

\begin{document}
\title[Hermete-Hadamard's type inequalities ]{On the generalization some
intgeral inequalities and their applications }
\author{Mehmet Zeki Sar\i kaya$^{\star }$}
\address{Department of Mathematics, Faculty of Science and Arts, D\"{u}zce
University, D\"{u}zce, Turkey}
\email{sarikayamz@gmail.com}
\thanks{$^{\star }$corresponding author}
\author{Nesip Aktan}
\email{nesipaktan@gmail.com}

\begin{abstract}
In this paper, a general integral identity for convex functions is derived.
Then, we establish new some inequalities of the Simpson and the
Hermite-Hadamard's type for functions whose absolute values of derivatives
are convex. Some applications for special means of real numbers are also
provided.
\end{abstract}

\maketitle

\section{Introduction}

Let $f:I\subseteq \mathbb{R\rightarrow R}$ be a convex mapping defined on
the interval $I$ of real numbers and $a,b\in I$ with $a<b.$ The following
double inequality:%
\begin{equation}
f(\frac{a+b}{2})\leq \frac{1}{b-a}\int\limits_{a}^{b}f(x)dx\leq \frac{%
f(a)+f(b)}{2}  \label{0}
\end{equation}%
is known in the literature as Hadamard inequality for convex mapping. Note
that some of the classical inequalities for means can be derived from (\ref%
{0}) for appropriate particular selections of the mapping $f.$ Both
inequalities hold in the reversed direction if $f$ is concave.

It is well known that the Hermite-Hadamard's inequality plays an important
role in nonlinear analysis. Over the last decade, this classical inequality
has been improved and generalized in a number of ways; there have been a
large number of research papers written on this subject, (see, \cite{SSDRPA}-%
\cite{Dragomir},\cite{Hussain}-\cite{K}, \cite{OAS}-\cite{sarikaya}, \cite%
{yang} ) and the references therein.

In \cite{sarikaya}, Sarikaya et. al. established inequalities for twice
differentiable convex mappings which are connected with Hadamard's
inequality, and they used the following lemma to prove their results:

\begin{lemma}
Let $f:I^{\circ }\subset \mathbb{R}\rightarrow \mathbb{R}$ be twice
differentiable function on $I^{\circ }$, $a,b\in I^{\circ }$($I^{\circ }$ is
the interior of $I$) with $a<b.$ If $f^{\prime \prime }\in L_{1}[a,b]$, then%
\begin{equation*}
\begin{array}{l}
\dfrac{1}{b-a}\dint_{a}^{b}f(x)dx-f(\dfrac{a+b}{2}) \\ 
\\ 
\ \ \ \ \ \ \ \ \ \ =\dfrac{\left( b-a\right) ^{2}}{2}\dint_{0}^{1}m\left(
t\right) \left[ f^{\prime \prime }(ta+(1-t)b)+f^{\prime \prime }(tb+(1-t)a)%
\right] dt,%
\end{array}%
\end{equation*}%
where%
\begin{equation*}
m(t):=\left\{ 
\begin{array}{ll}
t^{2} & ,t\in \lbrack 0,\frac{1}{2}) \\ 
&  \\ 
\left( 1-t\right) ^{2} & ,t\in \lbrack \frac{1}{2},1].%
\end{array}%
\right.
\end{equation*}
\end{lemma}

Also, the main inequalities in \cite{sarikaya}, pointed out as follows:

\begin{theorem}
Let $f:I\subset \mathbb{R}\rightarrow \mathbb{R}$ be twice differentiable
function on $I^{\circ }$ with $f^{\prime \prime }\in L_{1}[a,b]$. If $%
\left\vert f^{\prime \prime }\right\vert $ is a convex on $[a,b],$\ then%
\begin{equation}
\begin{array}{l}
\left\vert \dfrac{1}{b-a}\dint_{a}^{b}f(x)dx-f(\dfrac{a+b}{2})\right\vert
\leq \dfrac{\left( b-a\right) ^{2}}{24}\left[ \dfrac{\left\vert f^{\prime
\prime }\left( a\right) \right\vert +\left\vert f^{\prime \prime }\left(
b\right) \right\vert }{2}\right] .%
\end{array}
\label{H1}
\end{equation}
\end{theorem}

\begin{theorem}
Let $f:I\subset \mathbb{R}\rightarrow \mathbb{R}$ be twice differentiable
function on $I^{\circ }$ such that $f^{\prime \prime }\in L_{1}[a,b]$ where $%
a,b\in I,$ $a<b$. If $\left\vert f^{\prime \prime }\right\vert ^{q}$ is a
convex on $[a,b],$\ $q\geq 1$, then%
\begin{equation}
\begin{array}{l}
\begin{array}{l}
\left\vert \dfrac{1}{b-a}\dint_{a}^{b}f(x)dx-f(\dfrac{a+b}{2})\right\vert
\leq \dfrac{\left( b-a\right) ^{2}}{24}\left[ \dfrac{\left\vert f^{\prime
\prime }\left( a\right) \right\vert ^{q}+\left\vert f^{\prime \prime }\left(
b\right) \right\vert ^{q}}{2}\right] ^{\frac{1}{q}}.%
\end{array}%
\end{array}
\label{H2}
\end{equation}%
where $\frac{1}{p}+\frac{1}{q}=1.$
\end{theorem}

In \cite{sarikaya1}, Sarikaya et. al. prove some inequalities related to
Simpson's inequality for functions whose derivatives in absolute value at
certain powers are convex functions:

\begin{theorem}
\label{t.2.2} Let $f:I\subset \mathbb{R}\rightarrow \mathbb{R}$ be twice
differentiable mapping on $I^{\circ }$ such that $f^{\prime \prime }\in L_{1}%
\left[ a,b\right] ,\ $where $a,b\in I$ with $a<b.$ If $\left\vert f^{\prime
\prime }\right\vert ^{q}$ is a convex on $\left[ a,b\right] $ and $q\geq 1,$
then the following inequality holds:%
\begin{multline*}
\left\vert \dfrac{1}{6}\left[ f(a)+4f\left( \dfrac{a+b}{2}\right) +f(b)%
\right] -\dfrac{1}{b-a}\dint_{a}^{b}f(x)dx\right\vert \\
\\
\leq \left( b-a\right) ^{2}\left( \dfrac{1}{162}\right) ^{1-\frac{1}{q}}%
\left[ \left( \frac{59\left\vert f^{\prime \prime }\left( a\right)
\right\vert ^{q}+133\left\vert f^{\prime \prime }\left( b\right) \right\vert
^{q}}{3^{5}.2^{7}}\right) ^{\frac{1}{q}}+\left( \frac{133\left\vert
f^{\prime \prime }\left( a\right) \right\vert ^{q}+59\left\vert f^{\prime
\prime }\left( b\right) \right\vert ^{q}}{3^{5}.2^{7}}\right) ^{\frac{1}{q}}%
\right]
\end{multline*}%
where $\frac{1}{p}+\frac{1}{q}=1.$
\end{theorem}

In recent years many authors have studied error estimations for Simpson's
inequality; for refinements, counterparts, generaliations and new Simpson's
type inequalities, see \cite{ADD}, \cite{DAC}  \cite{sarikaya1} and \cite%
{sarikaya2}.

In this paper, in order to provide a unified approach to establish midpoint
inequality, trapezoid inequality and Simpson's inequality for functions
whose absolute values of derivatives are convex, we will derive a general
integral identity for convex functions. Finally, some applications for
special means of real numbers are provided.

\section{Main Results}

In order to prove our main theorems, we need the following Lemma:

\begin{lemma}
\label{z} Let $I\subset \mathbb{R}$ be an open interval, $a,b\in I$ with $%
a<b.$ If $f:I\rightarrow \mathbb{R}$ is a twice differentiable mapping such
that \ $f^{\prime \prime }$is integrable and $0\leq \lambda \leq 1.$ Then
the following identity holds:%
\begin{equation}
(\lambda -1)f(\frac{a+b}{2})-\lambda \frac{f(a)+f(b)}{2}+\dfrac{1}{b-a}%
\dint_{a}^{b}f(x)dx=\left( b-a\right) ^{2}\dint_{0}^{1}k(t)f^{\prime \prime
}(ta+(1-t)b)dt  \label{1z}
\end{equation}%
where 
\begin{equation*}
k(t)=\left\{ 
\begin{array}{ll}
\frac{1}{2}t(t-\lambda ), & 0\leq t\leq \frac{1}{2} \\ 
&  \\ 
\frac{1}{2}(1-t)(1-\lambda -t), & \frac{1}{2}\leq t\leq 1.%
\end{array}%
\right.
\end{equation*}
\end{lemma}

\begin{proof}
It suffices to note that%
\begin{eqnarray}
I &=&\dint_{0}^{1}k\left( t\right) f^{\prime \prime }(ta+(1-t)b)dt  \notag \\
&&  \notag \\
&=&\frac{1}{2}\dint_{0}^{\frac{1}{2}}t(t-\lambda )f^{\prime \prime
}(ta+(1-t)b)dt+\frac{1}{2}\dint_{\frac{1}{2}}^{1}(1-t)(1-\lambda
-t)f^{\prime \prime }(ta+(1-t)b)dt  \notag \\
&&  \notag \\
&=&I_{1}+I_{2}.  \label{2z}
\end{eqnarray}%
Integrating by parts twice, we can state:%
\begin{eqnarray}
I_{1} &=&\frac{1}{2}\dint_{0}^{\frac{1}{2}}t(t-\lambda )f^{\prime \prime
}(ta+(1-t)b)dt  \notag \\
&&  \notag \\
&=&\dfrac{1}{2}t(t-\lambda )\dfrac{f^{\prime }(ta+(1-t)b)}{a-b}\underset{0}{%
\overset{\frac{1}{2}}{\mid }}-\dfrac{1}{2\left( a-b\right) }\dint_{0}^{\frac{%
1}{2}}\left( 2t-\lambda \right) f^{\prime }(ta+(1-t)b)dt  \notag \\
&&  \label{3z} \\
&=&\dfrac{1}{4(b-a)}(\lambda -\frac{1}{2})f^{\prime }(\frac{a+b}{2})+\frac{%
(\lambda -1)}{2(b-a)^{2}}f(\frac{a+b}{2})  \notag \\
&&  \notag \\
&&-\frac{\lambda }{2(b-a)^{2}}f(b)+\dfrac{1}{\left( b-a\right) ^{2}}%
\dint_{0}^{\frac{1}{2}}f(ta+(1-t)b)dt,  \notag
\end{eqnarray}%
and similarly, we get%
\begin{eqnarray}
I_{2} &=&\frac{1}{2}\dint_{\frac{1}{2}}^{1}(1-t)(1-\lambda -t)f^{\prime
\prime }(ta+(1-t)b)dt  \notag \\
&&  \notag \\
&=&-\dfrac{1}{4(b-a)}(\lambda -\frac{1}{2})f^{\prime }(\frac{a+b}{2})+\frac{%
(\lambda -1)}{2(b-a)^{2}}f(\frac{a+b}{2})  \label{4z} \\
&&  \notag \\
&&-\frac{\lambda }{2(b-a)^{2}}f(a)+\dfrac{1}{\left( b-a\right) ^{2}}\dint_{%
\frac{1}{2}}^{1}f(ta+(1-t)b)dt.  \notag
\end{eqnarray}%
Using (\ref{3z}) and (\ref{4z}) in (\ref{2z}), it follows that%
\begin{equation*}
I=\dfrac{1}{\left( b-a\right) ^{2}}\left[ (\lambda -1)f(\frac{a+b}{2}%
)-\lambda \frac{f(a)+f(b)}{2}+\dfrac{1}{b-a}\dint_{0}^{1}f(ta+(1-t)b)dt%
\right] .
\end{equation*}%
Thus, using the change of the variable $x=ta+(1-t)b$ for $t\in \left[ 0,1%
\right] $ and by multiplying the both sides by $\left( b-a\right) ^{2},$ we
have the conclusion (\ref{1z}).
\end{proof}

Using this Lemma we can obtain the following general integral inequalities:

\begin{theorem}
\label{z2} Let $I\subset \mathbb{R}$ be an open interval, $a,b\in I$ with $%
a<b$ and $f:I\rightarrow \mathbb{R}$ be twice differentiable mapping such
that \ $f^{\prime \prime }$is integrable and $0\leq \lambda \leq 1.$ If $%
\left\vert f^{\prime \prime }\right\vert $ is a convex on $\left[ a,b\right]
,$ then the following inequalities hold:%
\begin{equation}
\begin{array}{l}
\left\vert (\lambda -1)f(\dfrac{a+b}{2})-\lambda \dfrac{f(a)+f(b)}{2}+\dfrac{%
1}{b-a}\dint_{a}^{b}f(x)dx\right\vert  \\ 
\\ 
\leq \left\{ 
\begin{array}{ll}
\dfrac{\left( b-a\right) ^{2}}{12}\left[ \left( \lambda ^{4}+\left(
1+\lambda \right) (1-\lambda )^{3}+\dfrac{5\lambda -3}{4}\right) \left\vert
f^{\prime \prime }\left( a\right) \right\vert \right.  & \text{for }0\leq
\lambda \leq \frac{1}{2} \\ 
&  \\ 
\left. +\left( \lambda ^{4}+\left( 2-\lambda \right) \lambda ^{3}+\dfrac{%
1-3\lambda }{4}\right) \left\vert f^{\prime \prime }\left( b\right)
\right\vert \right] , &  \\ 
\begin{array}{l}
\\ 
\end{array}
&  \\ 
\dfrac{\left( b-a\right) ^{2}\left( 3\lambda -1\right) }{48}\left[
\left\vert f^{\prime \prime }\left( a\right) \right\vert +\left\vert
f^{\prime \prime }\left( b\right) \right\vert \right]  & \text{for }\frac{1}{%
2}\leq \lambda \leq 1.%
\end{array}%
\right. 
\end{array}
\label{5z}
\end{equation}
\end{theorem}

\begin{proof}
From \ Lemma \ref{z} and by definition of $k(t),$ we get%
\begin{eqnarray}
&&\left\vert (\lambda -1)f(\frac{a+b}{2})-\lambda \frac{f(a)+f(b)}{2}+\dfrac{%
1}{b-a}\dint_{a}^{b}f(x)dx\right\vert   \notag \\
&&  \notag \\
&\leq &\left( b-a\right) ^{2}\dint_{0}^{1}\left\vert k(t)\right\vert
\left\vert f^{\prime \prime }(ta+(1-t)b)\right\vert dt  \label{6z} \\
&&  \notag \\
&=&\frac{\left( b-a\right) ^{2}}{2}\left\{ \dint_{0}^{\frac{1}{2}}\left\vert
t(t-\lambda )\right\vert \left\vert f^{\prime \prime }(ta+(1-t)b)\right\vert
dt+\dint_{\frac{1}{2}}^{1}\left\vert (1-t)(1-\lambda -t)\right\vert
\left\vert f^{\prime \prime }(ta+(1-t)b)\right\vert dt\right\}   \notag \\
&&  \notag \\
&=&\frac{\left( b-a\right) ^{2}}{2}\{J_{1}+J_{2}\}.  \notag
\end{eqnarray}

We assume that $0\leq \lambda \leq \frac{1}{2}$, then using the convexity of 
$\left\vert f^{\prime \prime }\right\vert ,$ we get%
\begin{eqnarray}
J_{1} &\leq &\dint_{0}^{\frac{1}{2}}\left\vert t(t-\lambda )\right\vert %
\left[ t\left\vert f^{\prime \prime }(a)\right\vert +(1-t)\left\vert
f^{\prime \prime }(b)\right\vert \right] dt  \notag \\
&&  \notag \\
&=&\dint_{0}^{\lambda }t(\lambda -t)\left[ t\left\vert f^{\prime \prime
}(a)\right\vert +(1-t)\left\vert f^{\prime \prime }(b)\right\vert \right]
dt+\dint_{\lambda }^{\frac{1}{2}}t(t-\lambda )\left[ t\left\vert f^{\prime
\prime }(a)\right\vert +(1-t)\left\vert f^{\prime \prime }(b)\right\vert %
\right] dt  \label{7z} \\
&&  \notag \\
&=&\left[ \frac{\lambda ^{4}}{6}+\frac{3-8\lambda }{3.2^{6}}\right]
\left\vert f^{\prime \prime }(a)\right\vert +\left[ \frac{\left( 2-\lambda
\right) \lambda ^{3}}{6}+\frac{5-16\lambda }{3.2^{6}}\right] \left\vert
f^{\prime \prime }(b)\right\vert ,  \notag
\end{eqnarray}%
and similarly, we have%
\begin{eqnarray}
J_{2} &\leq &\dint_{\frac{1}{2}}^{1-\lambda }(1-t)(1-\lambda -t)\left[
t\left\vert f^{\prime \prime }(a)\right\vert +(1-t)\left\vert f^{\prime
\prime }(b)\right\vert \right] dt  \notag \\
&&  \notag \\
&&+\dint_{1-\lambda }^{1}(1-t)(t+\lambda -1)\left[ t\left\vert f^{\prime
\prime }(a)\right\vert +(1-t)\left\vert f^{\prime \prime }(b)\right\vert %
\right] dt  \label{8z} \\
&&  \notag \\
&=&\left[ \frac{1+\lambda }{6}(1-\lambda )^{3}+\frac{48\lambda -27}{3.2^{6}}%
\right] \left\vert f^{\prime \prime }(a)\right\vert +\left[ \frac{\lambda
^{4}}{6}+\frac{3-8\lambda }{3.2^{6}}\right] \left\vert f^{\prime \prime
}(b)\right\vert .  \notag
\end{eqnarray}%
Using (\ref{7z}) and (\ref{8z}) in (\ref{6z}), we see that the first
inequality of (\ref{5z}) holds.

On the other hand, let $\frac{1}{2}\leq \lambda \leq 1,$ then, using the
convexity of $\left\vert f^{\prime \prime }\right\vert $ and by simple
computation we have%
\begin{eqnarray}
J_{1} &\leq &\dint_{0}^{\frac{1}{2}}\left\vert t(t-\lambda )\right\vert %
\left[ t\left\vert f^{\prime \prime }(a)\right\vert +(1-t)\left\vert
f^{\prime \prime }(b)\right\vert \right] dt  \notag \\
&&  \notag \\
&=&\dint_{0}^{\frac{1}{2}}t(\lambda -t)\left[ t\left\vert f^{\prime \prime
}(a)\right\vert +(1-t)\left\vert f^{\prime \prime }(b)\right\vert \right] dt
\label{2e} \\
&&  \notag \\
&=&\frac{8\lambda -3}{3.2^{6}}\left\vert f^{\prime \prime }(a)\right\vert +%
\frac{16\lambda -5}{3.2^{6}}\left\vert f^{\prime \prime }(b)\right\vert , 
\notag
\end{eqnarray}%
and similarly,%
\begin{eqnarray}
J_{2} &\leq &\dint_{\frac{1}{2}}^{1}\left\vert (1-t)(1-\lambda
-t)\right\vert \left[ t\left\vert f^{\prime \prime }(a)\right\vert
+(1-t)\left\vert f^{\prime \prime }(b)\right\vert \right] dt  \notag \\
&&  \notag \\
&=&\dint_{\frac{1}{2}}^{1}(1-t)(t+\lambda -1)\left[ t\left\vert f^{\prime
\prime }(a)\right\vert +(1-t)\left\vert f^{\prime \prime }(b)\right\vert %
\right] dt  \label{3e} \\
&&  \notag \\
&=&\frac{16\lambda -5}{3.2^{6}}\left\vert f^{\prime \prime }(a)\right\vert +%
\frac{8\lambda -3}{3.2^{6}}\left\vert f^{\prime \prime }(b)\right\vert . 
\notag
\end{eqnarray}%
Thus, if we use the (\ref{2e}) and (\ref{3e}) in (\ref{6z}), we obtain the
second inequality of (\ref{5z}). This completes the proof.
\end{proof}

Another similar result may be extended in the following theorem:

\begin{theorem}
\label{z3} Let $I\subset \mathbb{R}$ be an open interval, $a,b\in I$ with $%
a<b$ and $f:I\rightarrow \mathbb{R}$ be twice differentiable mapping such
that \ $f^{\prime \prime }$is integrable and $0\leq \lambda \leq 1.$ If $%
\left\vert f^{\prime \prime }\right\vert ^{q}$ is a convex on $\left[ a,b%
\right] ,\ q\geq 1,$ then the following inequalities hold:%
\begin{equation}
\begin{array}{l}
\left\vert (\lambda -1)f(\dfrac{a+b}{2})-\lambda \dfrac{f(a)+f(b)}{2}+\dfrac{%
1}{b-a}\dint_{a}^{b}f(x)dx\right\vert \\ 
\\ 
\leq \left\{ 
\begin{array}{ll}
\dfrac{\left( b-a\right) ^{2}}{2}\left( \dfrac{\lambda ^{3}}{3}+\dfrac{%
1-3\lambda }{24}\right) ^{1-\frac{1}{q}} & \ \text{for }0\leq \lambda \leq 
\frac{1}{2} \\ 
&  \\ 
\times \left\{ \left( \left[ \dfrac{\lambda ^{4}}{6}+\dfrac{3-8\lambda }{%
3.2^{6}}\right] \left\vert f^{\prime \prime }(a)\right\vert ^{q}+\left[ 
\dfrac{\left( 2-\lambda \right) \lambda ^{3}}{6}+\dfrac{5-16\lambda }{3.2^{6}%
}\right] \left\vert f^{\prime \prime }(b)\right\vert ^{q}\right) ^{\frac{1}{q%
}}\right. &  \\ 
&  \\ 
+\left. \left( \left[ \dfrac{1+\lambda }{6}(1-\lambda )^{3}+\dfrac{48\lambda
-27}{3.2^{6}}\right] \left\vert f^{\prime \prime }(a)\right\vert ^{q}+\left[ 
\dfrac{\lambda ^{4}}{6}+\dfrac{3-8\lambda }{3.2^{6}}\right] \left\vert
f^{\prime \prime }(b)\right\vert ^{q}\right) ^{\frac{1}{q}}\right\} , &  \\ 
\begin{array}{l}
\\ 
\end{array}
&  \\ 
\dfrac{\left( b-a\right) ^{2}}{2}\left( \dfrac{3\lambda -1}{24}\right) ^{1-%
\frac{1}{q}}\left\{ \left( \dfrac{8\lambda -3}{3.2^{6}}\left\vert f^{\prime
\prime }(a)\right\vert ^{q}+\dfrac{16\lambda -5}{3.2^{6}}\left\vert
f^{\prime \prime }(b)\right\vert ^{q}\right) ^{\frac{1}{q}}\right. & \text{%
for }\frac{1}{2}\leq \lambda \leq 1, \\ 
&  \\ 
\left. +\left( \dfrac{16\lambda -5}{3.2^{6}}\left\vert f^{\prime \prime
}(a)\right\vert ^{q}+\dfrac{8\lambda -3}{3.2^{6}}\left\vert f^{\prime \prime
}(b)\right\vert ^{q}\right) ^{\frac{1}{q}}\right\} & 
\end{array}%
\right.%
\end{array}
\label{9z}
\end{equation}%
\ where $\frac{1}{p}+\frac{1}{q}=1$.
\end{theorem}

\begin{proof}
Suppose that $q\geq 1.$ From Lemma \ref{z} and using the well known power
mean inequality, we have%
\begin{eqnarray}
&&\left\vert (\lambda -1)f(\frac{a+b}{2})-\lambda \frac{f(a)+f(b)}{2}+\dfrac{%
1}{b-a}\dint_{a}^{b}f(x)dx\right\vert  \notag \\
&&  \notag \\
&\leq &\left( b-a\right) ^{2}\dint_{0}^{1}\left\vert k\left( t\right)
\right\vert \left\vert f^{\prime \prime }\left( tb+\left( 1-t\right)
a\right) \right\vert dt  \notag \\
&&  \notag \\
&\leq &\frac{\left( b-a\right) ^{2}}{2}\left\{ \dint_{0}^{\frac{1}{2}%
}\left\vert t(t-\lambda )\right\vert \left\vert f^{\prime \prime
}(ta+(1-t)b)\right\vert dt+\dint_{\frac{1}{2}}^{1}\left\vert (1-t)(1-\lambda
-t)\right\vert \left\vert f^{\prime \prime }(ta+(1-t)b)\right\vert dt\right\}
\notag \\
&&  \label{10z} \\
&=&\frac{\left( b-a\right) ^{2}}{2}\left\{ \left( \dint_{0}^{\frac{1}{2}%
}\left\vert t(t-\lambda )\right\vert dt\right) ^{1-\frac{1}{q}}\left(
\dint_{0}^{\frac{1}{2}}\left\vert t(t-\lambda )\right\vert \left\vert
f^{\prime \prime }(ta+(1-t)b)\right\vert ^{q}dt\right) ^{\frac{1}{q}}\right.
\notag \\
&&  \notag \\
&&+\left. \left( \dint_{\frac{1}{2}}^{1}\left\vert (1-t)(1-\lambda
-t)\right\vert dt\right) ^{1-\frac{1}{q}}\left( \dint_{\frac{1}{2}%
}^{1}\left\vert (1-t)(1-\lambda -t)\right\vert \left\vert f^{\prime \prime
}(ta+(1-t)b)\right\vert ^{q}dt\right) ^{\frac{1}{q}}\right\} .  \notag
\end{eqnarray}%
Let $0\leq \lambda \leq \frac{1}{2}.$ Then, since $\left\vert f^{\prime
}\right\vert ^{q}$ is convex on $\left[ a,b\right] ,$ we know that for $t\in
\lbrack 0,1]$%
\begin{equation*}
\left\vert f^{\prime }(ta+(1-t)b)\right\vert ^{q}\leq t\left\vert f^{\prime
}(a)\right\vert ^{q}+(1-t)\left\vert f^{\prime }(b)\right\vert ^{q},
\end{equation*}%
hence, by simple computation%
\begin{eqnarray}
&&\dint_{0}^{\frac{1}{2}}\left\vert t(t-\lambda )\right\vert \left\vert
f^{\prime \prime }(ta+(1-t)b)\right\vert ^{q}dt  \notag \\
&&  \notag \\
&\leq &\dint_{0}^{\lambda }t(\lambda -t)\left[ t\left\vert f^{\prime \prime
}(a)\right\vert ^{q}+(1-t)\left\vert f^{\prime \prime }(b)\right\vert ^{q}%
\right] dt+\dint_{\lambda }^{\frac{1}{2}}t(t-\lambda )\left[ t\left\vert
f^{\prime \prime }(a)\right\vert ^{q}+(1-t)\left\vert f^{\prime \prime
}(b)\right\vert ^{q}\right] dt  \notag \\
&&  \label{11z} \\
&=&\left[ \frac{\lambda ^{4}}{6}+\frac{3-8\lambda }{3.2^{6}}\right]
\left\vert f^{\prime \prime }(a)\right\vert ^{q}+\left[ \frac{\left(
2-\lambda \right) \lambda ^{3}}{6}+\frac{5-16\lambda }{3.2^{6}}\right]
\left\vert f^{\prime \prime }(b)\right\vert ^{q},  \notag
\end{eqnarray}%
\begin{eqnarray}
&&\dint_{\frac{1}{2}}^{1}\left\vert (1-t)(1-\lambda -t)\right\vert
\left\vert f^{\prime \prime }(ta+(1-t)b)\right\vert ^{q}dt  \notag \\
&&  \notag \\
&\leq &\dint_{\frac{1}{2}}^{1-\lambda }(1-t)(1-\lambda -t)\left[ t\left\vert
f^{\prime \prime }(a)\right\vert ^{q}+(1-t)\left\vert f^{\prime \prime
}(b)\right\vert ^{q}\right] dt  \notag \\
&&  \label{12z} \\
&&+\dint_{1-\lambda }^{1}(1-t)(t+\lambda -1)\left[ t\left\vert f^{\prime
\prime }(a)\right\vert ^{q}+(1-t)\left\vert f^{\prime \prime }(b)\right\vert
^{q}\right] dt  \notag \\
&&  \notag \\
&=&\left[ \frac{1+\lambda }{6}(1-\lambda )^{3}+\frac{48\lambda -27}{3.2^{6}}%
\right] \left\vert f^{\prime \prime }(a)\right\vert ^{q}+\left[ \frac{%
\lambda ^{4}}{6}+\frac{3-8\lambda }{3.2^{6}}\right] \left\vert f^{\prime
\prime }(b)\right\vert ^{q},  \notag
\end{eqnarray}%
\begin{equation}
\dint_{0}^{\frac{1}{2}}\left\vert t(t-\lambda )\right\vert
dt=\dint_{0}^{\lambda }t(\lambda -t)dt+\dint_{\lambda }^{\frac{1}{2}%
}t(t-\lambda )dt=\frac{\lambda ^{3}}{3}+\frac{1-3\lambda }{24},  \label{13z}
\end{equation}%
and%
\begin{equation}
\dint_{\frac{1}{2}}^{1}\left\vert (1-t)(1-\lambda -t)\right\vert dt=\dint_{%
\frac{1}{2}}^{1-\lambda }(1-t)(1-\lambda -t)dt+\dint_{1-\lambda
}^{1}(1-t)(t+\lambda -1)dt=\frac{\lambda ^{3}}{3}+\frac{1-3\lambda }{24}.
\label{14z}
\end{equation}%
Thus, using (\ref{11z})-(\ref{14z}) in (\ref{10z}), we obtain the first
inequality of (\ref{9z}).

Now, let $\frac{1}{2}\leq \lambda \leq 1,$ then, using the convexity of $%
\left\vert f^{\prime \prime }\right\vert ^{q}$, we have%
\begin{eqnarray}
J_{1} &\leq &\dint_{0}^{\frac{1}{2}}\left\vert t(t-\lambda )\right\vert %
\left[ t\left\vert f^{\prime \prime }(a)\right\vert ^{q}+(1-t)\left\vert
f^{\prime \prime }(b)\right\vert ^{q}\right] dt  \notag \\
&&  \notag \\
&=&\dint_{0}^{\frac{1}{2}}t(\lambda -t)\left[ t\left\vert f^{\prime \prime
}(a)\right\vert ^{q}+(1-t)\left\vert f^{\prime \prime }(b)\right\vert ^{q}%
\right] dt  \label{5e} \\
&&  \notag \\
&=&\frac{8\lambda -3}{3.2^{6}}\left\vert f^{\prime \prime }(a)\right\vert
^{q}+\frac{16\lambda -5}{3.2^{6}}\left\vert f^{\prime \prime }(b)\right\vert
^{q},  \notag
\end{eqnarray}%
similarly,%
\begin{eqnarray}
J_{2} &\leq &\dint_{\frac{1}{2}}^{1}\left\vert (1-t)(1-\lambda
-t)\right\vert \left[ t\left\vert f^{\prime \prime }(a)\right\vert
^{q}+(1-t)\left\vert f^{\prime \prime }(b)\right\vert ^{q}\right] dt  \notag
\\
&&  \notag \\
&=&\dint_{\frac{1}{2}}^{1}(1-t)(t+\lambda -1)\left[ t\left\vert f^{\prime
\prime }(a)\right\vert ^{q}+(1-t)\left\vert f^{\prime \prime }(b)\right\vert
^{q}\right] dt  \label{6e} \\
&&  \notag \\
&=&\frac{16\lambda -5}{3.2^{6}}\left\vert f^{\prime \prime }(a)\right\vert
^{q}+\frac{8\lambda -3}{3.2^{6}}\left\vert f^{\prime \prime }(b)\right\vert
^{q},  \notag
\end{eqnarray}%
and so,%
\begin{equation}
\dint_{\frac{1}{2}}^{1}\left\vert (1-t)(1-\lambda -t)\right\vert
dt=\dint_{0}^{\frac{1}{2}}\left\vert t(t-\lambda )\right\vert dt=\dint_{0}^{%
\frac{1}{2}}t(\lambda -t)dt=\frac{3\lambda -1}{24}.  \label{7e}
\end{equation}%
Thus, if we use the (\ref{5e}), (\ref{6e}) and (\ref{7e}) in (\ref{10z}), we
obtain the second inequality of (\ref{9z}). This completes the proof.
\end{proof}

\section{Applications to Quadrature Formulas}

In this section we point out some particular inequalities which generalize
some classical results such as : trapezoid inequality, Simpson's inequality,
midpoint inequality and others.

\begin{proposition}[Midpoint inequality]
Under the assumptions Theorem \ref{z2} with $\lambda =0$ in Theorem \ref{z2}
, then we get the following inequality,%
\begin{equation*}
\left\vert \dfrac{1}{b-a}\dint_{a}^{b}f(x)dx-f(\frac{a+b}{2})\right\vert
\leq \dfrac{\left( b-a\right) ^{2}}{24}\left[ \frac{\left\vert f^{\prime
\prime }\left( a\right) \right\vert +\left\vert f^{\prime \prime }\left(
b\right) \right\vert }{2}\right] .
\end{equation*}
\end{proposition}

\begin{proposition}[Trapezoid inequality]
Under the assumptions Theorem \ref{z2} with $\lambda =1$ in Theorem \ref{z2}%
, then we have%
\begin{equation*}
\left\vert \dfrac{1}{b-a}\dint_{a}^{b}f(x)dx-\frac{f(a)+f(b)}{2}\right\vert
\leq \dfrac{\left( b-a\right) ^{2}}{12}\left[ \frac{\left\vert f^{\prime
\prime }\left( a\right) \right\vert +\left\vert f^{\prime \prime }\left(
b\right) \right\vert }{2}\right] .
\end{equation*}
\end{proposition}

\begin{proposition}[Simpson inequality]
\label{c1} Under the assumptions Theorem \ref{z2} with $\lambda =\frac{1}{3}$
in Theorem \ref{z2}, then we get%
\begin{equation*}
\left\vert \frac{1}{6}\left[ f(a)+4f(\frac{a+b}{2})+f(b)\right] -\dfrac{1}{%
b-a}\dint_{a}^{b}f(x)dx\right\vert \leq \dfrac{\left( b-a\right) ^{2}}{168}%
\left[ \left\vert f^{\prime \prime }\left( a\right) \right\vert +\left\vert
f^{\prime \prime }\left( b\right) \right\vert \right] .
\end{equation*}
\end{proposition}

\begin{proposition}
Under the assumptions Theorem \ref{z2} with $\lambda =\frac{1}{2}$ in
Theorem \ref{z2} , then we get%
\begin{equation*}
\left\vert \dfrac{1}{b-a}\dint_{a}^{b}f(x)dx-\frac{1}{2}\left[ f(\frac{a+b}{2%
})+\frac{f(a)+f(b)}{2}\right] \right\vert \leq \dfrac{\left( b-a\right) ^{2}%
}{48}\left[ \frac{\left\vert f^{\prime \prime }\left( a\right) \right\vert
+\left\vert f^{\prime \prime }\left( b\right) \right\vert }{2}\right] .
\end{equation*}
\end{proposition}

\begin{proposition}
\label{c2} Under assumptations Theorem \ref{z3} with $\lambda =0$ in Theorem %
\ref{z3} , then we get the following "midpoint inequality",%
\begin{multline*}
\left\vert \dfrac{1}{b-a}\dint_{a}^{b}f(x)dx-f(\frac{a+b}{2})\right\vert \\
\\
\leq \dfrac{\left( b-a\right) ^{2}}{48}\left[ \left( \frac{3\left\vert
f^{\prime \prime }\left( a\right) \right\vert ^{q}+5\left\vert f^{\prime
\prime }\left( b\right) \right\vert ^{q}}{8}\right) ^{\frac{1}{q}}+\left( 
\frac{5\left\vert f^{\prime \prime }\left( a\right) \right\vert
^{q}+3\left\vert f^{\prime \prime }\left( b\right) \right\vert ^{q}}{8}%
\right) ^{\frac{1}{q}}\right] .
\end{multline*}
\end{proposition}

\begin{proposition}
\label{c3} Under assumptations Theorem \ref{z3} with $\lambda =1$ in Theorem %
\ref{z3} , then we get"trapezoid inequality"%
\begin{multline*}
\left\vert \dfrac{1}{b-a}\dint_{a}^{b}f(x)dx-\frac{f(a)+f(b)}{2}\right\vert
\\
\\
\leq \dfrac{\left( b-a\right) ^{2}}{24}\left[ \left( \frac{5\left\vert
f^{\prime \prime }\left( a\right) \right\vert ^{q}+11\left\vert f^{\prime
\prime }\left( b\right) \right\vert ^{q}}{16}\right) ^{\frac{1}{q}}+\left( 
\frac{11\left\vert f^{\prime \prime }\left( a\right) \right\vert
^{q}+5\left\vert f^{\prime \prime }\left( b\right) \right\vert ^{q}}{16}%
\right) ^{\frac{1}{q}}\right] .
\end{multline*}
\end{proposition}

\begin{proposition}
\label{c4} Under assumptations Theorem \ref{z3} with $\lambda =\frac{1}{3}$
in Theorem \ref{z3} , then we get"Simpson inequality"%
\begin{multline*}
\left\vert \frac{1}{6}\left[ f(a)+4f(\frac{a+b}{2})+f(b)\right] -\dfrac{1}{%
b-a}\dint_{a}^{b}f(x)dx\right\vert \\
\\
\leq \dfrac{\left( b-a\right) ^{2}}{162}\left[ \left( \frac{59\left\vert
f^{\prime \prime }\left( a\right) \right\vert ^{q}+133\left\vert f^{\prime
\prime }\left( b\right) \right\vert ^{q}}{3.2^{6}}\right) ^{\frac{1}{q}%
}+\left( \frac{133\left\vert f^{\prime \prime }\left( a\right) \right\vert
^{q}+59\left\vert f^{\prime \prime }\left( b\right) \right\vert ^{q}}{3.2^{6}%
}\right) ^{\frac{1}{q}}\right] .
\end{multline*}
\end{proposition}

\section{Applications to Special Means}

We shall consider the following special means:

(a) The arithmetic mean: $A=A(a,b):=\dfrac{a+b}{2},$ \ $a,b\geq 0,$

(b) The geometric mean: $G=G(a,b):=\sqrt{ab},$ \ $a,b\geq 0,$

(c) The harmonic mean: 
\begin{equation*}
H=H\left( a,b\right) :=\dfrac{2ab}{a+b},\ a,b>0,
\end{equation*}

(d) The logarithmic mean: 
\begin{equation*}
L=L\left( a,b\right) :=\left\{ 
\begin{array}{ccc}
a & if & a=b \\ 
&  &  \\ 
\frac{b-a}{\ln b-\ln a} & if & a\neq b%
\end{array}%
\right. \text{, \ \ \ }a,b>0,
\end{equation*}

(e) The Identric mean:%
\begin{equation*}
I=I\left( a,b\right) :=\left\{ 
\begin{array}{ccc}
a & \text{if} & a=b \\ 
&  &  \\ 
\frac{1}{e}\left( \frac{b^{b}}{a^{a}}\right) ^{\frac{1}{b-a}}\text{ } & 
\text{if} & a\neq b%
\end{array}%
\right. \text{, \ \ \ }a,b>0,
\end{equation*}

(f) The $p-$logarithmic mean

\begin{equation*}
L_{p}=L_{p}(a,b):=\left\{ 
\begin{array}{ccc}
\left[ \frac{b^{p+1}-a^{p+1}}{\left( p+1\right) \left( b-a\right) }\right] ^{%
\frac{1}{p}} & \text{if} & a\neq b \\ 
&  &  \\ 
a & \text{if} & a=b%
\end{array}%
\right. \text{, \ \ \ }p\in \mathbb{R\diagdown }\left\{ -1,0\right\} ;\;a,b>0%
\text{.}
\end{equation*}%
It is well known \ that $L_{p}$ is monotonic nondecreasing \ over $p\in 
\mathbb{R}$ with $L_{-1}:=L$ and $L_{0}:=I.$ In particular, we have the
following inequalities%
\begin{equation*}
H\leq G\leq L\leq I\leq A.
\end{equation*}%
Now, using the results of Section 3, some new inequalities is derived for
the above means.

\begin{proposition}
\label{p.1} Let $a,b\in R$, $0<a<b$ and $n\in \mathbb{N}$, $n>2.$ Then, we
have%
\begin{equation*}
\left\vert \frac{1}{3}A\left( a^{n},b^{n}\right) +\frac{2}{3}A^{n}\left(
a,b\right) -L_{n}^{n}\left( a,b\right) \right\vert \leq n(n-1)\frac{\left(
b-a\right) ^{2}}{168}\left[ a^{n-2}+b^{n-2}\right] .
\end{equation*}
\end{proposition}

\begin{proof}
The assertion follows from Proposition \ref{c1} applied to convex mapping $%
f\left( x\right) =x^{n},\;x\in R.$
\end{proof}

\begin{proposition}
\label{p.2} Let $a,b\in R$, $0<a<b.$ Then, for all $q\geq 1$, we have%
\begin{multline*}
\left\vert L^{-1}\left( a,b\right) -A^{-1}\left( a,b\right) \right\vert \\
\leq n(n-1)\frac{\left( b-a\right) ^{2}}{48}\left\{ \left( \frac{%
3a^{(n-2)q}+5b^{(n-2)q}}{8}\right) ^{\frac{1}{q}}+\left( \frac{%
5a^{(n-2)q}+3b^{(n-2)q}}{8}\right) ^{\frac{1}{q}}\right\}
\end{multline*}%
and%
\begin{multline*}
\left\vert L^{-1}\left( a,b\right) -H^{-1}\left( a,b\right) \right\vert \\
\leq n(n-1)\frac{\left( b-a\right) ^{2}}{24}\left\{ \left( \frac{%
5a^{(n-2)q}+11b^{(n-2)q}}{16}\right) ^{\frac{1}{q}}+\left( \frac{%
11a^{(n-2)q}+5b^{(n-2)q}}{16}\right) ^{\frac{1}{q}}\right\} .
\end{multline*}
\end{proposition}

\begin{proof}
The assertion follows from Proposition \ref{c2} and Proposition \ref{c3}
applied to the convex mapping $f\left( x\right) =1/x,\;x\in \left[ a,b\right]
,$ respectively$.$
\end{proof}

\begin{proposition}
\label{p.3} Let $a,b\in R$, $0<a<b.$ Then, for all $q\geq 1$, we have%
\begin{eqnarray*}
&&\left\vert \frac{1}{3}H^{-1}\left( a,b\right) +\frac{2}{3}A^{-1}\left(
a,b\right) -L^{-1}\left( a,b\right) \right\vert \\
&\leq &\dfrac{\left( b-a\right) ^{2}}{162}\left\{ \left( \dfrac{59}{3.2^{6}}%
\left\vert \frac{2}{b^{3}}\right\vert ^{q}+\dfrac{133}{3.2^{6}}\left\vert 
\frac{2}{a^{3}}\right\vert ^{q}\right) ^{\frac{1}{q}}+\left( \dfrac{133}{%
3.2^{6}}\left\vert \frac{2}{b^{3}}\right\vert ^{q}+\dfrac{59}{3.2^{6}}%
\left\vert \frac{2}{a^{3}}\right\vert ^{q}\right) ^{\frac{1}{q}}\right\} .
\end{eqnarray*}
\end{proposition}

\begin{proof}
The assertion follows from Proposition \ref{c4} applied to the convex
mapping $f\left( x\right) =1/x,\;x\in \left[ a,b\right] .$
\end{proof}


\begin{thebibliography}{99}
\bibitem{ADD} M. Alomari, M. Darus and S.S. Dragomir, New inequalities of
Simpson's type for s-convex \ functions with applications, \textit{RGMIA
Res. Rep. Coll.}, 12 (4) (2009), Article 9.

\bibitem{SSDRPA} S.S. Dragomir and R.P. Agarwal, \textit{Two inequalities
for differentiable mappings and applications to special means of real
numbers and trapezoidal formula}, Appl. Math. Lett., 11(5) (1998), 91--95.

\bibitem{Drag1} S. S. Dragomir, \textit{Two mappings in connection to
Hadamard's inequalities}, J. Math. Anal. Appl. 167 (1992), 49--56.

\bibitem{Drag2} S. S. Dragomir, Y. J. Cho, and S. S. Kim, \textit{%
Inequalities of Hadamard's type for Lipschitzian mappings and their
applications}, J. Math. Anal. Appl. 245 (2000), 489--501.

\bibitem{Dragomir} S.S. Dragomir and C.E.M. Pearce, \textit{Selected Topics
on Hermite-Hadamard Inequalities and Applications}, RGMIA Monographs,
Victoria University, 2000. Online:[http://www.sta\textcurrency %
.vu.edu.au/RGMIA/monographs/hermite\_hadamard.html].

\bibitem{DAC} S.S. Dragomir, R.P. Agarwal and P. Cerone, On Simpson's
inequality and applications, \textit{J. of Inequal. Appl.}, 5(2000), 533-579.

\bibitem{Hussain} S. Hussain, M.I. Bhatti and M. Iqbal, \textit{%
Hadamard-type inequalities for }s\textit{-convex functions I,} Punjab Univ.
Jour. of Math\textit{.}, Vol.41, pp:51-60, (2009).

\bibitem{USK} U.S. K\i rmac\i , \textit{Inequalities for differentiable
mappings and applications to special means of real numbers and to midpoint
formula}, Appl. Math. Comp\textit{.,} 147 (2004), 137-146.

\bibitem{USKMEO} U.S. K\i rmac\i\ and M.E. \"{O}zdemir, \textit{On some
inequalities for differentiable mappings and applications to special means
of real numbers and to midpoint formula}, Appl. Math. Comp.\textit{,} 153
(2004), 361-368.

\bibitem{K} U.S. K\i rmac\i , \textit{Improvement and further generalization
of inequalities for differentiable mappings and applications}, Computers and
Math. with Appl\textit{.,} 55 (2008), 485-493.

\bibitem{LIU} B.Z. Liu, An inequality of Simpson type, \textit{Proc. R. Soc.
A, 461 (2005), 2155-2158.}

\bibitem{OAS} M. E. \"{O}zdemir, M. Avc\i\ and E. Set, \textit{On some
inequalities of Hermite--Hadamard type via m-convexity, }Appl. Math. Lett.
in press

\bibitem{CEMPJP} C.E.M. Pearce and J. Pe\v{c}ari\'{c}, \textit{Inequalities
for differentiable mappings with application to special means and quadrature
formulae}, Appl. Math. Lett., 13(2) (2000), 51--55.

\bibitem{sarikaya} M. Z. Sarikaya, A. Saglam and H. Y\i ld\i r\i m, \textit{%
New inequalities of Hermite-Hadamard type for functions whose second
derivatives absolute values are convex and quasi-convex}, arXiv:1005.0451,
submited.

\bibitem{sarikaya1} M. Z. Sarikaya, E. Set and M. E. \"{O}zdemir, \textit{On
New Inequalities of Simpson's Type for Functions whose Second Derivatives
Absolute Values are Convex, }RGMIA Res. Rep. Coll.,13(1) (2010), Supplement,
Article 1.

\bibitem{sarikaya2} M. Z. Sarikaya, E. Set and M. E. \"{O}zdemir, \textit{On
New Inequalities of Simpson's Type for Convex Functions, }RGMIA Res. Rep.
Coll.,13(2) (2010), Supplement, Article 2.

\bibitem{yang} G. S. Yang and K. L. Tseng, \textit{On certain integral
inequalities related to Hermite--Hadamard inequalities}, J. Math. Anal.
Appl. 239 (1999), 180--187.
\end{thebibliography}
\end{document}